\documentclass[a4paper]{amsart}

\author{Samuel Boissi\`{e}re} 

\title{On the McKay correspondences for the Hilbert scheme of points on the affine plane}

\address{Samuel Boissi\`{e}re, Fachbereich f\"{u}r Mathematik und Informatik, Staudinger Weg $9$,
Johannes Gutenberg-Universit\"{a}t Mainz, $55099$ Mainz, Deutschland}

\email{boissiere@mathematik.uni-mainz.de}

\urladdr{http://sokrates.mathematik.uni.mainz.de/$\sim$samuel}

\usepackage{amssymb}
\usepackage[all]{xy}
\usepackage[english]{babel}
\usepackage{youngtab}

\DeclareMathOperator{\Hilb}{Hilb} \DeclareMathOperator{\gr}{gr}
\DeclareMathOperator{\Coeff}{Coeff} \DeclareMathOperator{\Hom}{Hom}
\DeclareMathOperator{\codim}{codim} \DeclareMathOperator{\Supp}{Supp}
\DeclareMathOperator{\Spec}{Spec} \DeclareMathOperator{\Coh}{Coh}
\DeclareMathOperator{\rg}{rg} \DeclareMathOperator{\age}{age}

\newcommand{\IC}{\mathbb{C}}
\newcommand{\IN}{\mathbb{N}}
\newcommand{\IQ}{\mathbb{Q}}
\newcommand{\IZ}{\mathbb{Z}}

\newcommand{\hilb}{\Hilb^n(\IC^2)}
\newcommand{\Snhilb}{S_n\text{-}\Hilb\IC^{2n}}

\newcommand{\tH}{\widetilde{H}}
\newcommand{\tK}{\widetilde{K}}

\newcommand{\cC}{\mathcal{C}}
\newcommand{\cF}{\mathcal{F}}
\newcommand{\cO}{\mathcal{O}}

\theoremstyle{plain}
\newtheorem{theorem}{Theorem}[section]

\newtheorem{proposition}[theorem]{Proposition}
\newtheorem{lemma}[theorem]{Lemma}
\newtheorem{remark}[theorem]{Remark}
\newtheorem{corollary}[theorem]{Corollary}

\begin{document}

\begin{abstract}
The quotient of a finite-dimensional vector space by the action of a finite
subgroup of automorphisms is usually a singular variety. Under appropriate
assumptions, the McKay correspondence relates the geometry of nice resolutions
of singularities and the representations of the group. For the Hilbert scheme
of points on the affine plane, we study how different correspondences (McKay,
dual McKay and multiplicative McKay) are related to each other.
\end{abstract}



\maketitle

\pagestyle{myheadings}
\markboth{SAMUEL BOISSI\`{E}RE}{On the McKay correspondences for the Hilbert scheme of points on the affine plane}
%%%\markright{On the McKay correspondences for the Hilbert scheme of points on the affine plane}

\section{Introduction}

Let $V$ be a finite-dimensional complex vector space and $G\subset SL(V)$ a
finite subgroup of automorphisms. The quotient $V/G=\Spec \cO(V)^G$ is usually
a singular variety. \emph{McKay correspondences} aim to relate the geometry of
nice resolutions of singularities $Y\rightarrow V/G$ to the group $G$. Such
correspondences were first constructed by McKay \cite{McKay} and
Gonzales-Sprinberg-Verdier \cite{GSV} in dimension $2$ and then generalized in
dimension $3$ by several authors (for a survey of the subject, see Reid
\cite{Reid1} and references therein). In this paper, we are interested in
generalizations to higher dimensions dues to Bridgeland-King-Reid \cite{BKR},
Kaledin \cite{K3,K2,K1}, Ginzburg-Kaledin~\cite{GK} and Bezrukavnikov-Kaledin
\cite{BK}.

We suppose that the vector space $V$ is equipped with a symplectic form and
that the group $G$ preserves the symplectic structure. Suppose given a
\emph{crepant} resolution of singularities $Y\rightarrow V/G$ ($K_Y\cong
\cO_Y$). We study three constructions of McKay correspondences:

\noindent\textbf{McKay correspondence.} The derived categories of coherent
sheaves $D^b(\Coh(Y))$ and $D^b(\Coh_G(V))$ are equivalent
(Bezrukavnikov-Kaledin \cite{BK}). In some special cases, the \emph{Hilbert
scheme of regular $G$-orbits} $Y:=G\text{-}\Hilb V$ provides such a resolution
and the equivalence of categories is constructed as a Fourier-Mukai functor
(Bridgeland-King-Reid \cite{BKR}). This induces an isomorphism of Grothendieck
groups $K(Y)\cong K_G(V)$.

\noindent\textbf{Dual McKay correspondence.} Recall Kaledin's results
(\cite{K1}). The vector space $V$ is naturally stratified by the subspaces of
invariant vectors for different subgroups of $G$, inducing a stratification of
$V/G$. The resolution $Y\rightarrow V/G$ is \emph{semi-small} with respect to
this stratification. The maximal strata are indexed by conjugacy classes in $G$
and form a basis of the Borel-Moore homology with complex coefficients
$H^{\text{BM}}_*(Y)$. Denoting the space of $\IC$-valued functions on
$G$ invariant by conjugation by $\cC(G)$, we get a linear isomorphism
$H^{\text{BM}}_*(Y)\cong \cC(G)$.

\noindent\textbf{Multiplicative McKay correspondence.} The increasing
filtration of the group algebra $F^d\IC[G]:=\IC\{g\in G \,|\, \rg(id_V-g)\leq
d\}$ is compatible with the ring structure. By graduation and restriction to
the center $ZG$, we get a graded commutative algebra $\gr^FZG$. There is a
natural isomorphism of graded algebras $H^*(Y)\cong~ \gr^F ZG$
(Ginzburg-Kaledin \cite{GK}).

A natural question occurs at this point: how are these three correspondences
related to each other? In order to compare them, one is lead to the following
two problems, formulated by Ginzburg-Kaledin \cite{GK}:

\noindent\textbf{Poincar\'{e} duality problem.} (\cite[Problem $1.4$]{GK}) Compute
the Poincar\'{e} duality isomorphism:
$$
\gr^F ZG\longrightarrow H^*(Y) \overset{D}{\longrightarrow} H^{\text{BM}}_*(Y)
\longrightarrow \cC(G).
$$

\noindent\textbf{Chern character problem.} (\cite[Problem $1.5$]{GK}) Compute
the Chern character isomorphism:
$$
ZG\cong R(G)\underset{\IZ}{\otimes}\IC\cong K_G(V)\underset{\IZ}{\otimes}\IC
\longrightarrow K(Y)\underset{\IZ}{\otimes}\IC \overset{ch}{\longrightarrow}
H^*(Y) \longrightarrow \gr^F ZG.
$$

We study these questions in the particular case $V=\IC^n\otimes \IC^2$ with the
permutation action of the symmetric group $G=S_n$ and the canonical symplectic
structure. The Hilbert scheme $Y=\hilb$ provides a natural symplectic
resolution of singularities, isomorphic to the Hilbert scheme of regular orbits
$\Snhilb$. Here, the \emph{McKay correspondence} is realized by the
Bridgeland-King-Reid theorem \cite{BKR} and computed with Haiman's results
(\cite{H6,H5}); the multiplicative McKay correspondence is constructed by
Nakajima's operators (\cite{N3}) and the theorems of Lehn-Sorger \cite{LS} and
Vasserot \cite{V} and finally the \emph{dual McKay correspondence} is given by
natural subvarieties inducing a basis of the Borel-Moore homology. In this
context, we can replace the spaces $R(S_n)$, $\cC(S_n)$ and $ZS_n$ by the space
of \emph{symmetric functions} $\Lambda^n$ and work with rational coefficients.
We use the natural basis of \emph{Newton functions} $p_\lambda$ and \emph{Schur
functions} $s_\lambda$ indexed by \emph{partitions} $\lambda$ of $n$. We define
a graduation by $\deg p_\lambda:=n-l(\lambda)$ where $l(\lambda)$ is the
\emph{length} of the partition and consider the decreasing filtration
$F_d\Lambda^n:=\IQ\{p_\lambda\,|\,\deg p_\lambda \geq d\}$. For any symmetric
function $f\in \Lambda^n$ we denote the homogeneous component of
degree $k$ of $f$ by $[f]_k$.

In this setup, the \emph{Poincar\'{e} duality problem} consists in the computation
of the map:
$$
\gamma:\Lambda^n\longrightarrow H^*(\hilb,\IQ)
\overset{D}{\longrightarrow}H^{\text{BM}}_*(\hilb,\IQ) \longrightarrow
\Lambda^n
$$
and the \emph{Chern character problem} consists in the computation of the map:
$$
\Gamma:\Lambda^n\longrightarrow K(\hilb)\underset{\IZ}{\otimes}\IQ
\overset{ch}{\longrightarrow}H^*(\hilb,\IQ) \longrightarrow \Lambda^n.
$$

Denote the \emph{age} (or \emph{shifting degree}) of a partition
$\lambda=(\lambda_1,\lambda_2,\ldots)$ of $n$ by $\age(\lambda):=n-l(\lambda)$
and define his \emph{complexity degree} as $\langle \lambda \rangle
:=\prod\limits_{i\geq 1} \lambda_i$. These numbers can be given a more general
definition for any group acting on a vector space: the \emph{age} is computed
with a diagonalization of the action whereas the \emph{complexity degree} is
defined by the Frobenius decomposition of the action.

We prove the following formulas (\S \ref{ss:Poincare}, \S \ref{ss:Chern}):

\begin{theorem} For any partition $\lambda=(\lambda_1,\ldots,\lambda_k)$ of $n$,
\begin{align*}
\gamma(p_\lambda)&=\frac{(-1)^{\age(\lambda)}}{\langle \lambda \rangle }p_\lambda, \\
\Gamma (p_\lambda)&=\frac{(-1)^{\age(\lambda)}}{\langle \lambda \rangle}
p_\lambda +\sum_{\substack{\mu \vdash n
\\ l(\mu)<l(\lambda)}} g_{\lambda,\mu}p_\mu,
\end{align*}
for some coefficients $g_{\lambda,\mu}$.
\end{theorem}

\noindent Moreover we show the following theorem (\S \ref{ss:Conclusion}):

\begin{theorem}
The McKay correspondence is compatible with the decreasing topological
filtration of $K(\hilb)$ and the decreasing filtration of $\Lambda^n$.
\end{theorem}

Hence we can graduate the McKay correspondence and get a \emph{graded McKay
correspondence}. Our formulas then show that surprisingly, up to natural
identifications, the \emph{graded McKay correspondence}, the
\emph{multiplicative McKay correspondence} and the \emph{dual McKay
correspondence} are the same.

The main tool in our computations consists in computing combinatorial formulas
for the Chern classes of vector bundles on $\hilb$ linearizing the natural
action of the torus $\IC^*$ (\S \ref{ss:Chern classes}):

\begin{theorem}
Let $F$ be a $\IC^*$-linearized vector bundle of rank $r$ on $\hilb$ and
$f_1^\lambda,\ldots,f_r^\lambda$ the weights of the action on the fibre at each
fixed point. Then the Chern classes of $F$ written in $\Lambda^n$ are
$$
c_k(F)=\sum\limits_{\lambda\vdash n} \frac{1}{h(\lambda)}
\sigma_k(f_1^\lambda,\ldots,f_r^\lambda) [s_\lambda]_k,
$$
where the $\sigma_k(-)$ are the elementary symmetric functions.

The Chern characters of $F$ are
$$
ch_k(F)=\frac{1}{k!}\sum\limits_{\lambda\vdash n} \frac{1}{h(\lambda)}
\sum\limits_{i=1}^r \left(f_i^\lambda\right)^k [s_\lambda]_k.
$$
\end{theorem}

\section{Symmetric functions}
\label{s:Symmetric}

\subsection{The ring of symmetric functions} (\cite{McDo1,M})
\label{ss:Ring Symmetric}

Take independent indeterminates $x_1,\ldots,x_r$. Let the symmetric group $S_r$
act by permutation on $\IQ[x_1,\ldots,x_r]$ and denote the invariant ring by
$\Lambda_r:=~\IQ[x_1,\ldots,x_r]^{S_r}$. This ring is naturally graded by
degree and we denote the vector subspace generated by degree $n$
homogeneous symmetric polynomials by $\Lambda^n_r$. By adjoining other indeterminates, we may
construct the projective limit $\Lambda^n:=\underleftarrow{\lim}\Lambda^n_r$.
Then the \emph{ring of symmetric functions} is defined by
$\Lambda:=\bigoplus\limits_{n\geq 0}\Lambda^n$.

A \emph{partition} of an integer $n$ is a decreasing sequence of non-negative
integers $\lambda:=(\lambda_1,\ldots,\lambda_k)$ such that
$\sum\limits_{i=1}^k\lambda_i=n$ (we write $\lambda\vdash n$). The $\lambda_i$
are the \emph{parts} of the partition. If necessary, we extend a partition with
zero parts. The number $l(\lambda)$ of non-zero parts is the \emph{length} of
the partition and the sum $|\lambda|$ of the parts is the \emph{weight}. If a
partition $\lambda$ has $\alpha_1$ parts equal to $1$, $\alpha_2$ parts equal
to $2$, $\ldots$ we shall also denote it by
$\lambda:=(1^{\alpha_1},2^{\alpha_2},\ldots)$.

The \emph{Young diagram} of a partition $\lambda$ is defined by
$$
D(\lambda):=\{(i,j) \in \IN\times\IN \,|\, j<\lambda_{i+1} \}.
$$
In the representation of such a diagram, we follow a matrix convention:
$$
\Yvcentermath1\young(~xhh,~h~,~) 
\begin{array}{cc}
\lambda=(4,3,1) & x=(0,1)\\
|\lambda|=8 & h(x)=4 \\
l(\lambda)=3 &
\end{array}
$$

For each cell $x\in D(\lambda)$, the \emph{hook length} $h(x)$ at $x$ is the
number of cells on the right and below $x$. We shall also make use of the
number $n(\lambda):=\sum\limits_{i\geq 1} (i-1)\lambda_i$.

Set $p_0=1$ and define for $k\geq 1$ the \emph{power sum}
$p_k:=\sum\limits_{i\geq 1} x_i^k$. For a partition
$\lambda=(\lambda_1,\ldots,\lambda_k)$, the \emph{Newton function} is the
product $p_\lambda:=p_{\lambda_1}\cdots p_{\lambda_k}\in \Lambda^{|\lambda|}$.
The Newton functions form a basis of $\Lambda$ and
$\Lambda\cong\IQ[p_1,p_2,\ldots]$. Another natural basis of $\Lambda$ indexed
by partitions is given by the \emph{Schur functions} $s_\lambda$.

For a partition $\lambda=(1^{\alpha_1},2^{\alpha_2},\ldots)$, set
$z_\lambda:=\prod\limits_{r\geq 1} \alpha_r!r^{\alpha_r}$ and define a scalar
product on $\Lambda$ by $\langle p_\lambda,p_\mu\rangle=\delta_{\lambda,\mu}
z_\lambda$ where $\delta_{\lambda,\mu}$ is the Kronecker symbol. Then the basis
of Schur functions is orthonormal.

Let $\cC(S_n)$ be the $\IQ$-vector space of \emph{class functions} on $S_n$.
Since conjugacy classes in $S_n$ are indexed by partitions, the functions
$\chi_\lambda$ taking the value $1$ on the conjugacy class $\lambda$ and $0$
else form a basis of $\cC(S_n)$. Let $R(S_n)$ be the $\IQ$-vector space of
representations of $S_n$. By associating to each representation of $S_n$ his
\emph{character}, we get an isomorphism $\chi:R(S_n)\rightarrow \cC(S_n)$. The
\emph{Frobenius morphism} is the isomorphism $\Phi:\cC(S_n)\rightarrow
\Lambda^n$ characterized by $\Phi(\chi_\lambda)=z_\lambda^{-1}p_\lambda$. Let
$\chi^\lambda$ be the class function such that $\Phi(\chi^\lambda)=s_\lambda$
and $\chi^\lambda_\mu$ the value of $\chi^\lambda$ at the conjugacy class
$\mu$. The representations $V^\lambda$ of character $\chi^\lambda$ are the
irreducible representations of $S_n$ and we have the following base-change
formulas:
\begin{align*}
p_\mu&=\sum_{\lambda \vdash n} \chi^\lambda_\mu s_\lambda \text{ (Frobenius formula}), \\
s_\lambda&=\sum_{\mu \vdash n} z_\mu^{-1} \chi^\lambda_\mu p_\mu \text{
(inverse Frobenius formula)}.
\end{align*}

\subsection{Phethystic substitutions} (\cite{H1,McDo1})
\label{ss:Plethystic}

The identification $\Lambda=\IQ[p_1,p_2,\ldots]$ allows to specialize the
$p_k$'s to elements of any $\IQ$-algebra: the specialization extends uniquely
to an algebra homomorphism on $\Lambda$. For a formal Laurent series $E$ in
indeterminates $t_1,t_2,\ldots$ we define $p_k[E]$ to be the result of
replacing each indeterminate $t_i$ by $t_i^k$. Extending the specialization to
any symmetric function $f\in \Lambda$, we obtain the \emph{plethystic
substitution} of $E$ in $f$, denoted by $f[E]$. Our convention is that in a
plethystic substitution, $X$ stands for the sum of the original indeterminates
$x_1+x_2+\cdots$, so that $p_k[X]=p_k$.

\subsection{Macdonald polynomials}(\cite{H1,McDo2})
\label{ss:Macdonald}

We introduce indeterminates $q,t$ and consider the ring
$\Lambda_{\IQ(q,t)}:=\Lambda\underset{\IQ}{\otimes} \IQ(q,t)$. The scalar
product and the plethystic substitutions naturally extend to this situation.
For a partition $\mu$, we define $B_\mu(q,t):=\sum\limits_{(i,j)\in
D(\mu)}t^iq^j$. Set
$$
\Omega:=\exp\left(\sum\limits_{k\geq 1} \frac{p_k}{k}\right)
$$
and define a linear operator $\Delta:\Lambda_{\IQ(q,t)}\rightarrow
\Lambda_{\IQ(q,t)}$ by
$$
\Delta f=\left. f \left[ X+\frac{(1-t)(1-q)}{z} \right] \Omega[-zX]
\right|_{z^0}.
$$
The \emph{modified Macdonald polynomial} $\tH_\mu$ is the eigenvector of
$\Delta$ corresponding to the eigenvalue $1-(1-q)(1-t)B_\mu(q,t)$. These
polynomials form a basis of $\Lambda_{\IQ(q,t)}$ and decompose in the basis of
Schur functions as:
$$
\tH_\mu=\sum_{\lambda \vdash n} \tK_{\lambda,\mu}s_\lambda,
$$
where the $\tK_{\lambda,\mu}\in \IN[q,t]$ are called $q,t$-\emph{Kostka
polynomials}. We shall make use of the following specialization at $t=1/q$
(\cite[Proposition $3.5.10$]{H2}):
\begin{equation}\label{eq:evaluation McDo}
\tH_\mu(q,q^{-1})=q^{-n(\mu)}\prod_{x\in D(\mu)} \left(1-q^{h(x)}\right) s_\mu
\left[ \frac{X}{1-q}\right].
\end{equation}

\section{Hilbert schemes in the affine plane}
\label{s:Hilbert}

\subsection{Hilbert scheme of points}(\cite{B,EG,ES,Fo,G,LS,N1,N3,V})
\label{ss:Hilbert points}

The \emph{Hilbert scheme of $n$ points in the affine plane} $\hilb$ is the
smooth quasi-projective manifold of complex dimension $2n$ parameterizing
length $n$ finite subschemes in the plane $\IC^2$. The first projection
$B_n:=pr_{1*}\cO_{\Xi_n}$ of the universal family $\Xi_n\subset \hilb\times
\IC^2$ is the rank $n$ \emph{usual tautological bundle} on $\hilb$.

The manifold $\hilb$ has no odd singular cohomology; his even cohomology has no
torsion and is generated by algebraic cycles. We denote the
cohomology ring by $H^*(\hilb)$, the Borel-Moore homology by $H^{\text{BM}}_*(\hilb)$ and
the Grothendieck group of algebraic vector bundles (or equivalently of coherent
sheaves) by $K(\hilb)$, all with rational coefficients. 
Denote the Chern character by $ch:K(\hilb)\overset{\sim}{\longrightarrow} H^*(\hilb)$ and 
the Poincar\'{e} duality by $D:H^*(\hilb)\overset{\sim}{\longrightarrow} H^{\text{BM}}_*(\hilb)$.

Denote the Mumford quotient parameterizing length
$n$ effective zero-cycles in $\IC^2$ by $S^n\IC^2:=\IC^{2n}/S_n$. The \emph{Hilbert-Chow} morphism
$\rho:\hilb\rightarrow S^n\IC^2$ is a symplectic resolution of singularities,
semi-small with respect to the natural stratification $S_\lambda\IC^2:=\left\{
\sum\limits_{i=1}^k \lambda_i x_i \,|\, x_i\neq x_j\text{ for } i\neq
j\right\}$ for partitions $\lambda=(\lambda_1,\ldots,\lambda_k)$ of $n$. Each
subvariety $X_\lambda:=\rho^{-1}S_\lambda\IC^2$ is irreducible and locally
closed of dimension $n+l(\lambda)$.

There is a natural isomorphism
$$
\Psi:\Lambda^n \longrightarrow H^*(\hilb)
$$
constructed by use of geometric operators acting on the total sum of cohomology
of Hilbert schemes. For $i\geq 1$ denote by $X_{n,i}\subset \hilb\times
\Hilb^{n+i}(\IC^2)$ the subvariety of nested subschemes $\xi\subset\xi'$ such
that $\xi$ and $\xi'$ differ by a point of length $i$. Let $\pi_n,\pi_{n+i}$ be
the respective projections on $\hilb$ and $\Hilb^{n+i}(\IC^2)$ and define the
operator
$$
q_i:H^*(\hilb)\longrightarrow H^{*+2i-2}(\Hilb^{n+i}(\IC^2))
$$
by\footnote{For any continuous map $f:X\rightarrow Y$ between smooth oriented manifolds, we denote  the push-forward map in cohomology induced from the homological push-forward
by Poincar\'{e} duality by $f_!:H^*(X)\rightarrow H^*(Y)$.} 
$q_i(\alpha)=\pi_{n+i!}(\pi_n^*(\alpha)\cup[X_{n,i}])$. Denote the unit by
$|0\rangle\in H^0(\Hilb^{0}(\IC^2))$. Nakajima \cite{N3} proves that
the vectors
$$
q_\lambda:=q_{\lambda_1}\circ\cdots \circ q_{\lambda_k}|0\rangle\in
H^{2n-2l(\lambda)}(\hilb)
$$
where $\lambda=(\lambda_1,\ldots,\lambda_k)$ runs over all partitions of $n$
form a basis of $H^*(\hilb)$. The isomorphism $\Psi$ is defined by
$\Psi(p_\lambda)=q_\lambda$. We shall make use of the following observation
concerning the cohomology classes of the subvarieties $X_\lambda$:
\begin{equation}\label{eq:baseX}
\text{for } \lambda=(1^{\alpha_1},2^{\alpha_2},\ldots),\quad
\Big[\overline{X_\lambda} \Big]=\frac{1}{\prod\limits_{i\geq 1}\alpha_i!}
q_\lambda.
\end{equation}

We now describe the ring structure we have on $H^*(\hilb)$ as explained in
Lehn-Sorger \cite{LS} and Vasserot \cite{V}. Introduce the following graduation
of the group algebra $\IQ[S_n]$. For a permutation $\pi$ of cycle-type
$\lambda$, set $\deg(\pi):=n-l(\lambda)$ and denote the vector
subspace generated by degree $d$ elements by $\IQ[S_n](d)$. The natural ring structure on
$\IQ[S_n]$ is not compatible with the graduation but with the associated
increasing filtration:
$$
F^d\IQ[S_n]:=\bigoplus_{d'\leq d}\IQ[S_n](d').
$$
We consider the associated graded ring
$\gr^F\IQ[S_n]:=\bigoplus\limits_{d=0}^{n-1} F^d\IQ[S_n]/F^{d-1}\IQ[S_n]$. The
center $ZS_n$ of $\IQ[S_n]$ is generated by the homogeneous elements
$\chi_\lambda$ so it inherits graduation, filtration and ring structure. Denote by
$\gr^F\Lambda^n$ the corresponding ring via the Frobenius isomorphism. The
space $\gr^F\Lambda^n$ is graded by the \emph{cohomological degree} $\deg
p_\lambda:=n-l(\lambda)$ and the increasing filtration is denoted by
$$
F^d\Lambda^n:=\IQ\left\{p_\lambda\,|\,\deg p_\lambda\leq d\right\}.
$$
Then $\Psi:\gr^F\Lambda^n\longrightarrow H^{2*}(\hilb)$ is a isomorphism of
graded algebras (Lehn-Sorger \cite{LS}, Vasserot \cite{V}).

\subsection{Hilbert scheme of regular orbits} (\cite{BKR,H6,INj,IN,Nm})

The \emph{Hilbert scheme of $S_n$-regular orbits} $\Snhilb$ is defined as the
closure in the Hilbert scheme $\Hilb^{n!}(\IC^{2n})$ of the open set of
$S_n$-free orbits and is isomorphic to the Hilbert scheme $\hilb$ (Haiman
\cite[Theorem $5.1$]{H6}). Denote the universal family by $Z_n\subset
\Snhilb\times \IC^{2n}$ and set $P_n:=p_*\cO_{Z_n}$ considered as the rank $n!$
\emph{unusual tautological bundle} on $\hilb$. This bundle is equipped with a
natural $S_n$-action inducing the regular representation on each fiber.
Consider the diagram:
$$
\xymatrix{Z_n\ar[r]^-q \ar[d]_-p&\IC^{2n} \ar[d]\\
\hilb\ar[r]^-\rho&S^n\IC^2}
$$
Denote by $D^b(\hilb)$ the derived category of coherent sheaves and
$D^b_{S_n}(\IC^{2n})$ the derived category of coherent $S_n$-sheaves. In this
situation, we can apply the Bridgeland-King-Reid theorem (\cite{BKR}) and get
an isomorphism of Grothendieck groups
$$
\Upsilon:=q_!\circ p^!:K(\hilb)\rightarrow K_{S_n}(\IC^{2n}).
$$

Consider the following composition of vector space isomorphisms:
$$
\Theta:\Lambda^n\xrightarrow{\Phi^{-1}} \cC(S_n) \xrightarrow{\chi^{-1}} R(S_n)
\xrightarrow{\tau^{-1}} K_{S_n}(\IC^{2n})\xrightarrow{\Upsilon^{-1}} K(\hilb)
$$
where $\tau$ is the Thom isomorphism (here it is the restriction to a fibre,
see \cite[Theorem $5.4.17$]{CG}). The $S_n$-action on $P_n$ induces an
isotypical decomposition
$$
P_n=\bigoplus_{\mu\vdash n} \mathbf{V}^\mu\otimes P_\mu
$$
where $\mathbf{V}^\mu$ is the trivial bundle with fibre $V^\mu$ on $\hilb$ and
$P_\mu:=\Hom_{S_n}(\mathbf{V}^\mu,P_n)$. Then a easy computation similar to
\cite[Formula ($5.3$)]{IN} shows that:

\begin{proposition}\label{prop:IN} For $\mu\vdash n$, $\Theta(s_\mu)=P_\mu^*$.
In particular, the dual bundles $P_\mu^*$ form a basis of $K(\hilb)$.
\end{proposition}

\subsection{Torus action on the Hilbert scheme of points}(\cite{ES,H5})
\label{ss:Torus action}

The torus $T:=\IC^*$ acts on $\IC[x,y]$ by $s.x=sx, s.y=s^{-1}y$ for $s\in T$.
It induces a natural action on $\hilb$ with finitely many fixed points
$\xi_\lambda$ parameterized by the partitions $\lambda$ of $n$. The action
extends to all natural objets at issue over $\IC^2$.

Let $F$ be a $T$-linearized vector bundle on $\hilb$. Each fibre
$F(\xi_\lambda)$ has a structure of representation of $T$ and by identifying
the representation ring of $T$ with the ring of polynomials $R(T)\cong
\IZ[s,s^{-1}]$ we set
$$
F(\xi_\lambda):=F_\lambda(s):=\sum\limits_{i=1}^r s^{f_i^\lambda},
$$
where $f_i^\lambda\in \IZ$ are the \emph{weights} of the action of $T$ on
$F(\xi_\lambda)$.

In particular, we have the following result:

\begin{proposition}[Haiman]\emph{(\cite[Proposition $3.4$]{H5})}
\label{prop:fibreP} For $\lambda,\mu\vdash n$,
$$
P_\mu(\xi_\lambda)=\left.\tK_{\mu,\lambda}\right|_{t=s,q=s^{-1}}.
$$
\end{proposition}

\section{Chern classes of linearized bundles}
\label{s:Chern classes}

\subsection{Equivariant cohomology of the Hilbert scheme of points in the affine plane}(\cite{Br2,N2,V})
\label{ss:Equivariant cohomology}

Let $E_T\rightarrow B_T$ be the classifying bundle of $T$-vector bundles. For
any algebraic variety $X$ with an action of $T$, let $H_T^*(X)$ and $H^T_*(X)$
be the equivariant cohomology and the equivariant Borel-Moore homology with
rational coefficients. By definition, $H^*_T(X)=H^*(X_T)$ where $X_T:=(X\times
E_T)/B_T$. The ring $H^*_T(pt)$ is isomorphic to $\IQ[u]$ where $u$ is an
indeterminate of degree $2$ and $H^*_T(X)$ is a graded $\IQ[u]$-algebra. We
denote the unit by $1_X\in H^0_T(X)$ . If $X$ is smooth of pure dimension $d$, the
Poincar\'{e} duality $D:H^i_T(X)\rightarrow H^T_{d-i}(X)$ is an isomorphism and for
a proper $T$-equivariant morphism $f:Y\rightarrow X$ we have a push-forward morphism
$f_!:H^*_T(Y)\rightarrow H^*_T(X)$. In particular, any closed $T$-stable
subvariety $Y\overset{j}{\hookrightarrow} X$ defines a cohomology class
$[Y]_T:=j_!1_Y\in H^*_T(X)$. Any $T$-linearized vector bundle $F$ on $X$ has
$T$-equivariant Chern classes $c_k^T(F)$ and $T$-equivariant Chern characters
$ch_k^T(F)$ in $H^{2k}_T(X)$ such that if $j:X\hookrightarrow X_T$ is the
inclusion of a fibre, we have $c_k(F)=j^*c_k^T(F)$ and $ch_k(F)=j^*ch_k^T(F)$.
For any $\IQ[u]$-module $M$ we denote the localization of $M$ at the ideal
$\langle u-1\rangle$ by $M'$.

Let $\Sigma\subset \IC^2$ be the vertical axis and $X_{n,i}^\Sigma$ the
subvariety of $X_{n,i}$ whose points are nested subschemes $\xi\subset \xi'$
with extremal point on $\Sigma$. As in Vasserot \cite{V}, define for $i\geq 1$
the operator
$$
q_i^T[\Sigma]:H^*_T(\hilb)\longrightarrow H^{*+2i}_T(\Hilb^{n+i}(\IC^2))
$$
by $q_i^T[\Sigma](\alpha)=\pi_{n+i!}(\pi^*_n(\alpha) \cup [X_{n,i}^\Sigma]_T)$.
Vasserot \cite{V} proves that the vectors:
$$
q_\lambda^T[\Sigma]:=q_{\lambda_1}^T[\Sigma] \circ\cdots\circ
q_{\lambda_k}^T[\Sigma]1_{\Hilb^0(\IC^2)} \in H^{2n}_T(\hilb)
$$
for all partitions $\lambda=(\lambda_1,\ldots,\lambda_k)$ of $n$ form a basis
of $H^{2n}_T(\hilb)$ and constructs an isomorphism
$$
\phi:H^{2n}_T(\hilb)\longrightarrow \Lambda^n
$$
by\footnote{There is an inaccuracy in \cite{V} : a factor $z_{(i)}$ is missing
in the formula $(2)$ (see Nakajima \cite[Lemma~ $9.4$]{N1}).}
$\phi(q_\lambda^T[\Sigma])=p_\lambda$.

Define for $i\geq 1$ the operator
$$
q_i^T:H^*_T(\hilb)\longrightarrow H^{*+2i-2}_T(\Hilb^{n+i}(\IC^2))
$$
by $q_i^T(\alpha)=\pi_{n+i!}(\pi^*_n(\alpha)\cup[X_{n,i}]_T)$ and set
$$
q_\lambda^T:=q_{\lambda_1}^T\circ\cdots\circ q_{\lambda_k}^T 1_{\Hilb^0(\IC^2)}
\in H^{2n-2l(\lambda)}_T(\hilb).
$$
Since $[\Sigma]_T=u.[\IC^2]_T$ we see that
$q_\lambda^T[\Sigma]=u^{l(\lambda)}q_\lambda^T$. The inclusion of a fibre
$j:~\hilb\hookrightarrow(\hilb)_T$ gives $j^*q_\lambda^T=~q_\lambda$. We can
apply the Leray-Hirsch theorem to the situation:
$$
\hilb\overset{j}{\longrightarrow}(\hilb)_T\overset{p}{\longrightarrow}B_T
$$
and get an isomorphism of graded $H^*(B_T)$-moduls:
$$
H^*_T(\hilb)\cong H^*(B_T)\otimes H^*(\hilb).
$$
Since $B_T$ and $\hilb$ have no odd cohomology, the theorem gives a basis in
each cohomological degree:
$$
\begin{array}{ccc}
H^{2m}_T(\hilb)& \cong &\bigoplus\limits_{k=0}^{m} H^{2m-2k}(B_T)\otimes H^{2k}(\hilb)\\
u^{m-k} q_\lambda^T&\leftrightarrow & u^{m-k}\otimes q_\lambda
\end{array}
$$
where $\lambda$ is a partition of $n$ such that $n-l(\lambda)=k$.

The multiplication by $u$ sending $H^k_T(\hilb)$ to $H^{k+2}_T(\hilb)$ is
always injective. Since $H^q(\hilb)=0$ for $q\geq 2n$, the vector space
$H^{2n}_T(\hilb)$ contains all the information about the equivariant cohomology
and multiplication by $u$ becomes an isomorphism after this degree:
$$
H^0_T(\hilb) \overset{u\cdot}{\longrightarrow} \cdots
\overset{u\cdot}{\longrightarrow} H^{2n}_T(\hilb)
\underset{\sim}{\overset{u\cdot}{\longrightarrow}} H^{2n+2}_T(\hilb)
\underset{\sim}{\overset{u\cdot}{\longrightarrow}} \cdots .
$$
The Leray-Hirsch decomposition makes $H^{2n}_T(\hilb)$ a graded vector space:
we denote by $\gr H^{2n}_T(\hilb)$ the vector space with his graded structure. The
vectors $q_\lambda^T[\Sigma]$ form a homogeneous basis with $\deg
q_\lambda^T[\Sigma]=n-l(\lambda)$. By the choice of this basis, we have a
canonical isomorphism
$$
H^{2n}_T(\hilb)\xrightarrow{can.}\gr H^{2n}_T(\hilb)
$$
cutting up a vector in homogeneous components: for $\alpha\in H^{2n}_T(\hilb)$
we denote the component of degree $k$ in $\alpha$ by $\gr_k\alpha$. We have an
isomorphism of graded vector spaces
$$
J:\gr H^{2n}_T(\hilb) \longrightarrow H^*(\hilb)
$$
defined by $J(q_\lambda^T[\Sigma])=q_\lambda$. The morphism
$\phi:H^{2n}_T(\hilb)\rightarrow \Lambda^n$ also induces a isomorphism of
graded vector spaces $\gr\phi:\gr H^{2n}_T(\hilb)\rightarrow \gr\Lambda^n$ and
the following diagram is commutative:
$$
\xymatrix{\Lambda^n \ar[r]^-{can.} & \gr \Lambda^n \ar[dr]^-\Psi & \\
H^{2n}_T(\hilb) \ar[r] \ar[u]^-\phi \ar[r]^-{can.} & \gr H^{2n}_T(\hilb)
\ar[u]^-{\gr \phi} \ar[r]^-J & H^*(\hilb) }
$$
Define $[\lambda]\in H^{2n}_T(\hilb)$ by
$u^n\cdot[\lambda]=(-1)^nh(\lambda)^{-1}[\xi_\lambda]_T$ where $h(\lambda)$ is
the product of the hook lengths in the Young diagram of $\lambda$. We have:

\begin{proposition}[Vasserot]\emph{(\cite{V})}
\label{prop:Vasserot} For $\lambda\vdash n$, $\phi([\lambda])=s_\lambda$.
\end{proposition}

\subsection{Chern classes of linearized bundles} \text{}
\label{ss:Chern classes}

In the study of the Chern classes of natural vector bundles on $\hilb$, we
prove the following formulas:

\begin{theorem} \label{th:Chern classes}
Let $F$ be a $T$-linearized vector bundle of rank $r$ on $\hilb$ and
$f_1^\lambda,\ldots,f_r^\lambda$ the weights of the action on the fibre at each
fixed point. Then the Chern classes of $F$ written in $\Lambda^n$ via $\Psi$
are
$$
c_k(F)=\sum\limits_{\lambda\vdash n} \frac{1}{h(\lambda)}
\sigma_k(f_1^\lambda,\ldots,f_r^\lambda) \sum\limits_{\substack{\mu\vdash n \\
l(\mu)=n-k}} z_\mu^{-1}\chi^\lambda_\mu p_\mu,
$$
where the $\sigma_k(-)$ are the elementary symmetric functions.

The Chern characters of $F$ are
$$
ch_k(F)=\frac{1}{k!}\sum\limits_{\lambda\vdash n} \frac{1}{h(\lambda)}
\sum\limits_{i=1}^r \left(f_i^\lambda\right)^k \sum\limits_{\substack{\mu\vdash n \\
l(\mu)=n-k}} z_\mu^{-1}\chi^\lambda_\mu p_\mu.
$$
\end{theorem}

\begin{proof}
The inclusion of a fixed point is denoted by
$i_\lambda:\xi_\lambda\hookrightarrow\hilb$ and we set
$[\xi_\lambda]_T:=i_{\lambda!}1_{\xi_\lambda} \in H^{4n}_T(\hilb)$. The
inclusion of the fixed points locus is denoted by
$$
i_n:=\bigoplus\limits_{\lambda\vdash n} i_\lambda:(\hilb)^T\hookrightarrow
\hilb.
$$
By the localization theorem in equivariant cohomology, the direct image
$$
i_{n!}:H^*_T((\hilb)^T)'\rightarrow H^*_T(\hilb)'
$$
is an isomorphism. The inverse is given by
$$
\alpha\mapsto\sum\limits_{\lambda\vdash
n}\frac{i^*_\lambda\alpha}{c^T_{max}(T_{\xi_\lambda}\hilb)}1_{\xi_\lambda}.
$$
Let $\theta$ be the representation of $T$ of weight $1$. The isomorphism
$H^*_T(pt)~\cong \IQ[u]$ is given by the first Chern class so for $a\in \IZ$ we
have
$$
c_{tot}^T(\theta^{\otimes a})=1+auZ\in H^*(B_T)[Z],
$$
where $c_{tot}^T:=1+c_1^TZ+c_2^TZ^2+\cdots$ is the total Chern class. Then, by
the properties of the Chern classes we get:
$$
c_{tot}^T(i_\lambda^*F)=\prod\limits_{i=1}^r \left(1+f_i^\lambda uZ\right).
$$
In particular, we know (see Nakajima \cite{N1}) that the representation of the
fibre at $\xi_\lambda$ of the tangent space of $\hilb$  is given by
$$
T_{\xi_\lambda}\hilb\cong \bigoplus\limits_{x\in D(\lambda)}
\left(\theta^{h(x)}\oplus\theta^{-h(x)}\right).
$$
It follows that $ c_{max}^T(T_{\xi_\lambda}\hilb)=(-1)^nh(\lambda)^2u^{2n}$.
The inverse localization formula gives then in $H^*_T(\hilb)'[Z]$:
$$
c_{tot}^T(F)=(-1)^n\frac{1}{u^{2n}}\sum_{\lambda\vdash
n}\frac{1}{h(\lambda)^2}\prod\limits_{i=1}^r \left(1+f_i^\lambda
uZ\right)[\xi_\lambda]_T.
$$
Since $H^q(\hilb)=0$ for $q>2n$, we suppose $k\leq n$. The $Z^k$-term gives in
$H^*_T(\hilb)'$:
$$
u^{2n}c_k^T(F)=(-1)^nu^k\sum\limits_{\lambda\vdash n} \frac{1}{h(\lambda)^2}
\sigma_k(f_1^\lambda,\ldots,f_r^\lambda)[\xi_\lambda]_T,
$$
where the $\sigma_k(-)$ are the elementary symmetric functions. Since $u$ is
invertible in the localized module, we get in $H^{2n}_T(\hilb)'$:
$$
u^{n-k}c_k^T(F)=\sum\limits_{\lambda\vdash n} \frac{1}{h(\lambda)}
\sigma_k\left(f_1^\lambda,\ldots,f_r^\lambda\right)[\lambda].
$$
Since multiplication by $u$ is an isomorphism after $H^{2n}_T(\hilb)$, this
equation is in fact an equation in $H^{2n}_T(\hilb)$.

\begin{lemma}
$J\left(\gr_k\left(u^{n-k}c_k^T(F)\right)\right)=c_k(F)$.
\end{lemma}

\begin{proof}[Proof of the lemma]
With the Leray-Hirsch decomposition:
$$
H^{2k}_T(\hilb)\cong \bigoplus\limits_{j=0}^k H^{2j}(B_T)\otimes
H^{2k-2j}(\hilb)
$$
we can write $c_k^T(F)=\sum\limits_{j=0}^k u^j\otimes \alpha_j$ and
$u^{n-k}c_k^T(F)=\sum\limits_{j=0}^k u^{n-k+j}\otimes \alpha_j$. Then
$\gr_k\left(u^{n-k}c_k^T(F)\right)=u^{n-k}\otimes \alpha_0$ and
$J\left(\gr_k\left(u^{n-k}c_k^T(F)\right)\right)=\alpha_0$. Since $j^*u=0$, we
also have $c_k(F)=j^*c_k^T(F)=\alpha_0$.
\end{proof}
From this lemma, the commutativity of the diagram and the proposition
\ref{prop:Vasserot} we obtain the expression of the Chern classes of $F$ in
$\Lambda^n$:
$$
c_k(F)=\sum_{\lambda\vdash n} \frac{1}{h(\lambda)}
\sigma_k(f_1^\lambda,\ldots,f_r^\lambda)[s_\lambda]_k,
$$
where $[s_\lambda]_k$ means that we keep only the component of cohomological
degree $k$. The inverse Frobenius formula gives then:
$$
c_k(F)=\sum\limits_{\lambda\vdash n} \frac{1}{h(\lambda)}
\sigma_k(f_1^\lambda,\ldots,f_r^\lambda) \sum\limits_{\substack{\mu\vdash n \\
l(\mu)=n-k}} z_\mu^{-1}\chi^\lambda_\mu p_\mu.
$$
Similarly, starting from the formula (see also \cite{LQW})
$$
ch^T_k(i^*_\lambda F)=\frac{1}{k!}\sum_{i=1}^r (f_i^\lambda)^k,
$$
we find
$$
u^{n-k}ch_k^T(F)=\frac{1}{k!}\sum_{\lambda\vdash n}
\frac{1}{h(\lambda)}\sum_{i=1}^r \left(f_i^\lambda\right)^k [\lambda]
$$
and the naturality $j^*ch^T_k(F)=ch_k(F)$ implies in a similar manner:
$$
ch_k(F)=\frac{1}{k!}\sum_{\lambda\vdash n} \frac{1}{h(\lambda)}\sum_{i=1}^r
\left(f_i^\lambda\right)^k [s_\lambda]_k.
$$
\end{proof}

\begin{remark} The same method gives the Todd classes and with
these formulas we can recover the well-known formula \cite[Proposition
$5.2$]{LS} for the Chern classes of the tautological bundle $B_n$ over $\hilb$
(see \cite{Boissiere}).
\end{remark}

Let $t$ be an indeterminate, define $\omega_tp_k=t^{k-1}p_k$ and extend the
definition to an algebra homomorphism $\omega_t:\Lambda\rightarrow \Lambda[t]$.
Then $\omega_t p_\lambda=t^{|\lambda|-l(\lambda)}p_\lambda$: this notation
takes care of the cohomological degree. In particular,
$$
\omega_t s_\lambda=\sum_{k\geq 0} [s_\lambda]_k t^k.
$$

\begin{corollary} \label{co:Chern character}
The total Chern character of a $T$-linearized vector bundle $F$ on $\hilb$ is
$$
ch(F)=\sum_{\lambda\vdash n}
\frac{1}{h(\lambda)}\Coeff\left(t^0,\omega_ts_\lambda
F_\lambda(e^{1/t})\right).
$$
The total Chern character of the dual bundle $F^*$ is
$$
ch(F^*)=\sum_{\lambda\vdash n}
\frac{1}{h(\lambda)}\Coeff\left(t^0,\omega_ts_\lambda
F_\lambda(e^{-1/t})\right).
$$
\end{corollary}

\begin{proof}
By theorem \ref{th:Chern classes} and his proof, the total Chern character of
$F$ is:
$$
ch (F)=\sum_{\lambda \vdash n} \frac{1}{h(\lambda)} \sum_{i=1}^r \sum_{k\geq 0}
\frac{1}{k!} \left(f_i^\lambda\right)^k [s_\lambda]_k .
$$
By $F_\lambda(e^{1/t})=\sum\limits_{i=1}^r
e^{f_i^\lambda/t}=\sum\limits_{i=1}^r \sum\limits_{k\geq 0} \frac{1}{k!}
\left(f_i^\lambda\right)^k t^{-k}$ and $\omega_t s_\lambda = \sum\limits_{k\geq
0} [s_\lambda]_k t^k$ we deduce the first formula. The second formula is
similar since $ch_k (F^*)=(-1)^k ch_k (F)$.
\end{proof}

\section{Comparison problems}
\label{s:Comparison}

\subsection{Poincar\'{e} duality problem}\text{}
\label{ss:Poincare}

Denote by $\vartheta_{\overline{X_\mu}}\in H^{\text{BM}}_*(\hilb)$ the homology
fundamental class of the closed subvariety $\overline{X_\mu}$. By definition,
$D\Big[\overline{X_\mu}\Big]=\vartheta_{\overline{X_\mu}}$ and these classes
form the natural basis in homology as in Kaledin \cite{K1}. We get the
\emph{dual McKay correspondence} by defining a bijection
$H^{\text{BM}}_*(\hilb) \rightarrow \cC(S_n)$ with
$\Big[\overline{X_\mu}\Big]\mapsto \chi_\mu$. Composing with the Frobenius
morphism and introducing a sign (for a reason that will appear later) we
define:
$$
\begin{array}{cccc}\phi:&\Lambda^n&\rightarrow & H^{\text{BM}}_*(\hilb) \\
&p_\mu&\mapsto&(-1)^{n-l(\mu)}z_\mu \Big[\overline{X_\mu}\Big]
\end{array}
$$
With these notations, the `` Poincar\'{e} duality problem '' consists in the
computation of the dotted arrow $\gamma$:
$$
\xymatrix{ \Lambda^n \ar[rr]^-\phi && H^{\text{BM}}_*(\hilb) \\
\Lambda^n \ar@{.>}[u]^\gamma \ar[rr]^-\Psi &&  H^*(\hilb) \ar[u]^{D} }
$$

\begin{proposition}
\label{prop:Poincare problem} For $\mu\vdash n$,
$$
\gamma (p_\mu)=\frac{(-1)^{n-l(\mu)}}{\prod\limits_{i\geq 1}\mu_i} p_\mu.
$$
\end{proposition}

\begin{proof} Let
$\mu=(\mu_1,\ldots,\mu_k)=(1^{\alpha_1},2^{\alpha_2},\ldots)$ be a partition of
$n$. By definition, $\Psi(p_\mu)=q_\mu$ and with formula (\ref{eq:baseX}) we
have $\Psi(p_\mu)=\left(\prod\limits_{i\geq 1}\alpha_i!\right)
\Big[\overline{X_\mu}\Big]$. Since we have $z_\mu=\left(\prod\limits_{i\geq
1}\alpha_i!\right)\left(\prod\limits_{i\geq 1}\mu_i\right)$, we get the result.
\end{proof}

\subsection{Chern character problem}\text{}
\label{ss:Chern}

With our notations, the `` Chern character problem '' consists in the
computation of the dotted arrow $\Gamma$:
$$
\xymatrix{ \Lambda^n \ar[rr]^-\Psi && H^*(\hilb) \\
\Lambda^n \ar@{.>}[u]^\Gamma \ar[rr]^-\Theta &&  K(\hilb) \ar[u]^{ch} }
$$

\begin{theorem}
\label{th:Chern problem} For $\mu\vdash n$,
$$
\Gamma (p_\mu)=\frac{(-1)^{n-l(\mu)}}{\prod\limits_{i\geq 1}\mu_i} p_\mu
+\sum_{\substack{\nu\vdash n
\\ l(\nu)<l(\mu)}} g_{\mu,\nu}p_\nu
$$
for some coefficients $g_{\mu,\nu}$.
\end{theorem}

\begin{proof}
By the proposition \ref{prop:IN}, the map $\Gamma$ is characterized by:
$$
\Gamma(s_\mu)=ch (P^*_\mu).
$$
We use the proposition \ref{prop:fibreP} and we can apply corollary
\ref{co:Chern character} with $F=P_\mu$ and $F_\lambda(s)=\left.
\tK_{\mu,\lambda}\right|_{t=s,q=s^{-1}}$:
$$
\Gamma(s_\mu)=ch (P^*_\mu) = \sum_{\lambda \vdash n} \frac{1}{h(\lambda)}
\Coeff\left(t^0,\omega_t s_\lambda \tK_{\mu,\lambda}(e^{1/t},e^{-1/t})\right).
$$
Since $\tK_{\mu,\lambda}=\langle \tH_\lambda,s_\mu \rangle$ we get:
$$
\Gamma(s_\mu)= \sum_{\lambda \vdash n} \frac{1}{h(\lambda)} \Coeff\left(
t^0,\omega_t s_\lambda \langle \tH_\lambda(e^{1/t},e^{-1/t}),s_\mu \rangle
\right)
$$
and by linearity:
$$
\Gamma(p_\mu)= \sum_{\lambda \vdash n} \frac{1}{h(\lambda)}
\Coeff\left(t^0,\omega_t s_\lambda \langle \tH_\lambda(e^{1/t},e^{-1/t}),p_\mu
\rangle\right).
$$
We use formula (\ref{eq:evaluation McDo}):
$$
\tH_\lambda(q,q^{-1})=q^{-n(\lambda)}\prod_{x\in D(\lambda)}
\left(1-q^{h(x)}\right) s_\lambda \left[ \frac{X}{1-q}\right],
$$
and find:
$$
\textstyle{\Gamma(p_\mu)= \sum\limits_{\lambda \vdash n} \frac{1}{h(\lambda)}
\Coeff\left( t^0,\left.(\omega_t s_\lambda) q^{-n(\lambda)}\prod\limits_{x\in
D(\lambda)} \left(1-q^{h(x)}\right)\left\langle s_\lambda\left[
\frac{X}{1-q}\right],p_\mu \right\rangle \right|_{q=e^{1/t}} \right)}.
$$
Observe that:
$$
\left\langle s_\lambda \left[ \frac{X}{1-q} \right],p_\mu \right\rangle=
\left\langle s_\lambda,p_\mu \left[ \frac{X}{1-q}\right] \right\rangle  =
\prod_{i=1}^{l(\mu)}\frac{1}{(1-q^{\mu_i})} \chi^\lambda_\mu,
$$
so that in fact:
$$
\Gamma(p_\mu)= \sum_{\lambda \vdash n} \frac{1}{h(\lambda)} \Coeff\left(
t^0,(\omega_t s_\lambda) \chi^\lambda_\mu \left.\left( q^{-n(\lambda)}
\frac{\prod\limits_{x\in D(\lambda)}
(1-q^{h(x)})}{\prod\limits_{i=1}^{l(\mu)}(1-q^{\mu_i})}
\right)\right|_{q=e^{1/t}} \right).
$$
Since by construction:
$$
\omega_t s_\lambda=\sum_{\nu \vdash n} z_\nu^{-1} \chi^\lambda_\nu p_\nu
t^{n-l(\nu)},
$$
the decomposition of $\Gamma(p_\mu)$ in the basis $\{p_\nu\}$ is (with
$u=1/t$):
$$
\textstyle{\Gamma(p_\mu)=\sum\limits_{\nu \vdash n} z_\nu^{-1}
\sum\limits_{\lambda \vdash n} \frac{1}{h(\lambda)} \chi^\lambda_\nu
\chi^\lambda_\mu \Coeff \left.\left( u^{n-l(\nu)},\left( q^{-n(\lambda)}
\frac{\prod\limits_{x\in D(\lambda)}
(1-q^{h(x)})}{\prod\limits_{i=1}^{l(\mu)}(1-q^{\mu_i})} \right)\right|_{q=e^u}
\right)p_\nu}.
$$

\begin{lemma}
$$
\left.\left( q^{-n(\lambda)} \frac{\prod\limits_{x\in D(\lambda)}
(1-q^{h(x)})}{\prod\limits_{i=1}^{l(\mu)}(1-q^{\mu_i})} \right)\right|_{q=e^u}
= \frac{(-1)^{n-l(\mu)}h(\lambda)}{\prod\limits_{i=1}^{l(\mu)}\mu_i}
u^{n-l(\mu)} + \text{ upper powers }.
$$
\end{lemma}

\begin{proof}[Proof of the lemma]
By Taylor expansion we have:
\begin{align*}
1-e^{\mu_i u}&=-\mu_i u (1+u(\ldots)) ,\\
1-e^{h(x)u}&=-h(x)u(1+u(\ldots)),
\end{align*}
so the first term in the Taylor expansion of the expression is:
$$
\frac{\prod\limits_{x\in
D(\lambda)}(-h(x)u)}{\prod\limits_{i=1}^{l(\mu)}(-\mu_i
u)}=\frac{(-1)^{n-l(\mu)}h(\lambda)}{\prod\limits_{i=1}^{l(\mu)}\mu_i}
u^{n-l(\mu)}.
$$
\end{proof}

We deduce that if $l(\nu) > l(\mu)$ then the coefficient of $\Gamma(p_\mu)$ at
$p_\nu$ is zero. If $l(\nu)=l(\mu)$ then the coefficient of $\Gamma(p_\mu)$ at
$p_\nu$ is:
$$
g_{\mu,\nu}:=\frac{(-1)^{n-l(\mu)}}{\prod\limits_{i=1}^{l(\mu)}\mu_i}
z_\nu^{-1} \sum_{\lambda \vdash n} \chi^\lambda_\nu \chi^\lambda_\mu.
$$
In the Frobenius formula $p_\mu=\sum\limits_{\lambda \vdash n} \chi^\lambda_\mu
s_\lambda$, we take the scalar product with $p_\nu$ to find the identity
$$
\delta_{\mu,\nu}z_\nu=\sum_{\lambda \vdash n} \chi^\lambda_\nu
\chi^\lambda_\mu,
$$
which shows that $g_{\mu,\nu}=0$ if $\mu \neq \nu$ and
$g_{\mu,\mu}=\displaystyle{\frac{(-1)^{n-l(\mu)}}{\prod\limits_{i=1}^{l(\mu)}\mu_i}}$.
\end{proof}

\subsection{Conclusion}\text{}
\label{ss:Conclusion}

$\bullet$ The space $K(\hilb)$ has a decreasing \emph{topological} filtration
defined by the codimension of the support of coherent sheaves:
$$
F_dK(\hilb):=\IQ\left\{\cF\,|\,\codim \Supp \cF\geq d\right\}.
$$
The Chern character $ch:K(\hilb)\rightarrow H^*(\hilb)$ is compatible with this
filtration and the induced graded map $\gr ch:\gr K(\hilb)\rightarrow
H^*(\hilb)$ is given by the cohomology class of the support (see \cite[\S
$5.9$]{CG}):
$$
\gr ch(\cF)=[\Supp \cF].
$$
Denote the decreasing cohomological filtration on $\Lambda$ by:
$$
F_d\Lambda^n:=\IQ\left\{p_\lambda\,|\,\deg p_\lambda\geq d\right\},
$$
and $\gr_F\Lambda^n$ the induced graded vector space. Then our preceding
results mean:

\begin{theorem}
\label{th:filtration} The McKay correspondence $\Theta$ is compatible with the
topological filtration of $K(\hilb)$ and the decreasing filtration of
$\Lambda^n$.
\end{theorem}

\begin{proof}
By formula (\ref{eq:baseX}) the cohomology classes $\Big[\overline{X_\mu}
\Big]$ form a homogeneous basis of $H^*(\hilb)$, so the classes of the
structural sheaves $\cO_{\overline{X_\mu}}$ form a graded basis of $\gr
K(\hilb)$ with $\deg \cO_{\overline{X_\mu}}=n-l(\mu)$. Then the theorem
\ref{th:Chern problem} implies that
$\Theta^{-1}\left(\cO_{\overline{X_\mu}}\right)\in F_{n-l(\mu)}\Lambda^n$ and
after inversion of the matrix we get the result.
\end{proof}

$\bullet$ By theorem \ref{th:Chern problem}, the induced graded map:
$$
\gr\Gamma=\Psi\circ\gr ch \circ \gr\Theta:\gr_F\Lambda^n\rightarrow \gr
K(\hilb)\rightarrow H^*(\hilb)\rightarrow \gr^F\Lambda^n
$$
is defined by
$\gr\Gamma(p_\mu)=\displaystyle\frac{(-1)^{n-l(\mu)}}{\prod\limits_{i\geq
1}\mu_i} p_\mu$, which is exactly the same formula as for the map $\gamma$
(this justifies our sign modification in \S \ref{ss:Poincare}). The number
$\age(\lambda):=n-l(\lambda)$ is the \emph{age} (or \emph{shifting degree}) of
any $\sigma\in S_n$ of cycle-type $\lambda$ for the permutation action on
$\IC^{2n}$, defined as follows: the eigenvalues of $\sigma$ on $\IC^{2n}$ are
squares of the unit $e^{2\mathbf{i}\pi r_j}$, $r_j\in [0,1[$ and by definition
$\age(\sigma):=\sum\limits_{j\geq 1} r_j=n-l(\lambda)$. The number $\langle
\sigma \rangle:=\prod\limits_{i\geq 1}\lambda_i$ can be interpreted as the
product of the dimensions of the orbit vector subspaces one gets by performing
a Frobenius decomposition of the matrix of $\sigma$ on $\IC^{2n}$. We call this
number $\langle \sigma \rangle$ the \emph{complexity degree} of $\sigma$ (see
\cite{Boissiere}). Then, the map $p_\lambda\mapsto
\frac{(-1)^{\age(\lambda)}}{\langle \lambda \rangle}p_\lambda$ has a
signification for any finite group $G$ acting on a vector space $V$ : denoting
by $\chi_{[g]}$ the class function of the conjugacy class $[g]$ of $g\in G$, this
map $\cC(G)\rightarrow \cC(G)$ is $\chi_{[g]}\mapsto
\frac{(-1)^{\age(g)}}{\langle g\rangle}\chi_{[g]}$.

\nocite{*}
\bibliographystyle{amsplain}
\bibliography{BiblioMcKayHilbert}

\end{document}